\documentclass{article}

\usepackage{amsmath,amssymb,amscd}

\newtheorem{Thm}{Theorem}
\newtheorem{Lem}{Lemma}
\newtheorem{Prop}{Proposition}
\newtheorem{Rem}{Remark}

\newtheorem{Cor}{Corollary}

\begin{document}
\title{Applications of square roots of diffeomorphisms}
\author{Yoshihiro Sugimoto}
\date{}
\maketitle

\begin{abstract}
In this paper, we prove that on any contact manifold ${(M,\xi)}$, there exists an arbitrary ${C^{\infty}}$-small contactomorphism which does not admit a square root. In particular, there exists an arbitrary ${C^{\infty}}$-small contactomorphism which is not "autonomous". This paper is the first step to study the topology of ${Cont_0(M,\xi)\backslash \textrm{Aut}(M,\xi)}$. As an application, we also prove a similar result for the diffeomorphism group ${\textrm{Diff}(M)}$ for any smooth manifold $M$.
\end{abstract}

\section{Introduction}
For any closed manifold $M$, the set of diffeomorphisms ${\textrm{Diff}(M)}$ forms a group and any one-parameter subgroup ${f:\mathbb{R}\rightarrow \textrm{Diff}(M)}$ can be written in the following form.
\begin{equation*}
f(t)=\exp (tX)
\end{equation*}
Here, ${X\in \Gamma(TM)}$ is a vector field and ${\exp:\Gamma(TM)\rightarrow \textrm{Diff}(M)}$ is a time $1$ flow of a vector field. From the inverse function theorem, one might expect that there exists an open neighborhood of the zero section ${\mathcal{U}\subset \Gamma(TM)}$ such that 
\begin{equation*}
\exp:\mathcal{U}\longrightarrow \textrm{Diff}(M)
\end{equation*}
is a diffeomorphism onto an open neighborhood of ${\textrm{Id}\in \textrm{Diff}(M)}$. However, this is far from true (\cite{M}, Warning 1.6). So one might expect that the set of "autonomous" diffeomorphisms 
\begin{equation*}
\textrm{Aut}(M)=\exp(\Gamma(TM))
\end{equation*}
is a small subset of ${\textrm{Diff}(M)}$.

For a symplectic manifold ${(M,\omega)}$, the set of Hamiltonian diffeomorphisms ${\textrm{Ham}^c(M,\omega)}$ contains "autonomous" subset ${\textrm{Aut}(M,\omega)}$ as follows.
\begin{equation*}
\textrm{Aut}(M,\omega)=\biggl\{\exp(X) \ \biggl| \ \begin{matrix}X\ \textrm{is \ a \ time-independent \ Hamiltonian \ vector \ field} \\ \textrm{whose \ support\ is \ compact} \end{matrix} \biggl\}
\end{equation*}
In \cite{AF}, Albers and Frauenfelder proved that on any symplectic manifold there exists an arbitrary ${C^{\infty}}$-small Hamiltonian diffeomorphism not admitting a square root. In particular, there exists an arbitrary ${C^{\infty}}$-small Hamiltonian diffeomorphism in ${\textrm{Ham}^c(M,\omega)\backslash \textrm{Aut}(M,\omega)}$.
 
Polterovich and Shelukhin used spectral spread of Floer homology and Conley conjecture to prove that ${\textrm{Ham}^c(M,\omega)\backslash \textrm{Aut}(M,\omega)\subset \textrm{Ham}^c(M,\omega)}$ is ${C^{\infty}}$-dense and dense in the topology induced from Hofer's metric if ${(M,\omega)}$ is closed symplectically aspherical manifold (\cite{PS}). The author generalized this theorem to arbitrary closed symplectic manifolds and convex symplectic manifolds (\cite{S}).

One might expect that "contact manifold" version of these theorems hold. In this paper, we prove that there exists an arbitrary ${C^{\infty}}$-small contactomorphism not admitting a square root. In particular, there exists an arbitrary ${C^{\infty}}$-small contactomorphism in ${\textrm{Cont}^c_0(M,\xi)\backslash \textrm{Aut}(M,\xi)}$. So, this paper is a contact manifold version of \cite{AF}. As an application, we prove that there exists an arbitrary ${C^{\infty}}$-small diffeomorphism in ${\textrm{Diff}^c_0(M)}$ not admitting a square root. This also implies that  there exists an arbitrary ${C^{\infty}}$-small diffeomorphism in ${\textrm{Diff}^c_0(M)\backslash \textrm{Aut}(M)}$.

\section*{Acknowledgement}
The author thanks Kaoru Ono and Urs Frauenfelder for many useful comments, discussions and encouragement.

\section{Main result}
Let $M$ be a smooth ${(2n+1)}$-dimensional manifold without boundary. A $1$-form $\alpha$ on $M$ is called contact if ${(\alpha \wedge (d\alpha)^n)(p)\neq 0}$ holds on any ${p\in M}$. A codimension $1$ tangent distribution $\xi$ on $M$ is called contact structure if it is locally defined by ${\ker(\alpha)}$ for some (locally defined) contact form $\alpha$. A diffeomorphism ${\phi \in \textrm{Diff}(M)}$ is called contactomorphism if ${\phi_*\xi=\xi}$ holds (i.e. ${\phi}$ preserves the contact structure $\xi$). Let ${\textrm{Cont}^c_0(M,\xi)}$ be the set of compact supported contactomorphisms which are isotopic to ${\textrm{Id}}$. In other words, ${\textrm{Cont}^c_0(M,\xi)}$ is a connected component of compact supported contactomorphisms(${\textrm{Cont}^c(M,\xi)}$) which contains ${\textrm{Id}}$.

\begin{equation*}
Cont^c_0(M,\xi)=\biggl\{\phi_1 \ \biggl| \ \begin{matrix}\phi_t \ (t\in [0,1]) \ \textrm{is \ an \ isotopy \ of \ contactomorphisms}\\ \phi_0=\textrm{Id}, \ \cup_{t\in[0,1]}\textrm{supp}(\phi_t) \ \textrm{is \ compact}  \end{matrix}\biggl\}
\end{equation*}

Let ${X\in \Gamma^c(TM)}$ be a compact supported vector field on $M$. $X$ is called contact vector field if the flow of $X$ preserves the contact structure $\xi$ (i.e. ${\exp(X)_*\xi=\xi}$ holds). Let ${\Gamma^c_{\xi}(TM)}$ be the set of compact supported contact vector fields on $M$ and let ${\textrm{Aut}(M,\xi)}$ be their images.
\begin{equation*}
\textrm{Aut}(M,\xi)=\{\exp(X) \ | \ X\in \Gamma^c_{\xi}(TM)\}
\end{equation*}
We prove the following theorem.

\begin{Thm}
Let ${(M,\xi)}$ be a contact manifold without boundary. Let $\mathcal{W}$ be any ${C^{\infty}}$-open neighborhood of $\textrm{Id}\in \textrm{Cont}^c_0(M,\xi)$. Then, there exists ${\phi \in \mathcal{W}}$ such that
\begin{equation*}
\phi \neq \psi^2
\end{equation*}
holds for any ${\psi \in \textrm{Cont}^c_0(M,\xi)}$. In particular, ${\mathcal{W}\backslash \textrm{Aut}(M,\xi)}$ is not empty.
\end{Thm}

\begin{Rem}
If $\phi$ is autonomous ($\phi=\exp(X)$), $\phi$ has a square root ${\psi=\exp(\frac{1}{2}X)}$.
\end{Rem}

\begin{Cor}
The exponential map ${\exp:\Gamma^c_{\xi}(TM)\rightarrow \textrm{Cont}^c_{\xi}(M,\xi)}$ is not surjective.
\end{Cor}

We also consider diffeomorphism version of this theorem and corollary. Let $M$ be a smooth manifold without boundary and let ${\textrm{Diff}^c(M)}$ be the set of compact supported diffeomorhisms which contains ${\textrm{Id}}$.
\begin{equation*}
\textrm{Diff}^c(M)=\{\phi \in \textrm{Diff}(M) \ | \ \textrm{supp}(\phi) \ \textrm{is \ compact}\}
\end{equation*}

Let ${\textrm{Diff}^c_0(M)}$ be the connected component of ${\textrm{Diff}^c(M)}$(i.e. any element of ${\textrm{Diff}^c_0(M)}$ is isotopic to $\textrm{Id}$.). We define the set of autonomous diffeomorphisms as follows.
\begin{equation*}
\textrm{Aut}(M)=\{\exp (X) \ | \ X\in \Gamma^c(TM)\}
\end{equation*}
By combining the arguments in this paper and in \cite{AF}, we can prove the following theorem.

\begin{Thm}
Let $M$ be a smooth manifold without boundary. Let $\mathcal{W}$ be any ${C^{\infty}}$-open neighborhood of $\textrm{Id}\in \textrm{Diff}^c_0(M)$. Then, there exists ${\phi \in \mathcal{W} }$ such that  
\begin{equation*}
\phi\neq \psi^2
\end{equation*}
holds for any ${\psi \in \textrm{Diff}^c(M)}$. In particular, ${\mathcal{W}\backslash \textrm{Aut}(M)}$ is not empty.
\end{Thm} 

\begin{Cor}
The exponential map ${\exp:\Gamma^c(TM)\rightarrow \textrm{Diff}^c_0(M)}$ is not surjective.
\end{Cor}

\section{Milnor's criterion}
In this section, we review Milnor's arguments in \cite{M}. In \cite{M}, Milnor gave a criterion for the existence of square root of a diffeomorphism. We use this criterion later. We fix ${l\in \mathbb{N}_{\ge2}}$ and a diffeomorphism ${\phi\in \textrm{Diff}(M)}$. Let ${P^l(\phi)}$ be the set of "$l$-periodic orbits" as follows.
\begin{equation*}
P^l(\phi)=\{(x_1,\cdots,x_l) \ | \ x_i\neq x_j (i\neq j), x_j=\phi^{j-1}(x_1), x_1=\phi(x_l)\}/\sim
\end{equation*}
This equivalence relation $\sim$ is given by the natural ${\mathbb{Z}/l\mathbb{Z}}$-action
\begin{equation*}
(x_1,\cdots,x_l)\rightarrow (x_l,x_1,\cdots,x_{l-1}) .
\end{equation*}

\begin{Prop}[Milnor\cite{M}, Albers-Frauenfelder\cite{AF}]
Assume that ${\phi\in \textrm{Diff}(M)}$ has a square root (i.e. there exists ${\psi\in \textrm{Diff}(M)}$ such that ${\phi=\psi^2}$ holds). Then, there exists a free ${\mathbb{Z}/2\mathbb{Z}}$-action on ${P^{2k}(\phi)}$ (${k\in\mathbb{N}}$). In particular, ${\sharp P^{2k}(\phi)}$ is even if ${\sharp P^{2k}(\phi)}$ is finite.
\end{Prop}

The proof of this proposition is very simple. We explain the proof for the sake of self-containedness. Let ${\psi}$ be a square root of ${\phi}$ (${\phi=\psi^2}$). It is sufficient to prove that the natural ${\mathbb{Z}/2\mathbb{Z}}$-action
\begin{equation*}
[x_1,\cdots,x_{2k}]\longrightarrow [\psi(x_1),\cdots,\psi(x_{2k})]
\end{equation*} 
on ${P^{2k}(\phi)}$ is free. We prove this by contradiction. So, we assume that ${[x_1,\cdots,x_{2k}]}$ is a fixed point of this ${\mathbb{Z}/2\mathbb{Z}}$-action. Then there exists ${1\le r\le 2k}$ so that ${\psi(x_1)=x_r}$ holds. This implies that
\begin{gather*}
\psi(x_1)=x_r=\phi^{(r-1)}(x_1)=\psi^{2(r-1)}(x_1)  \\
x_1=\psi^{2r-3}(x_1)=\psi^{2r-3}(\psi^{2r-3}(x_1))=\phi^{2r-3}(x_1)
\end{gather*}
holds. Let $1\le l <2k$  be the greatest common divisor of ${2k}$ and ${2r-3}$. Then, 
\begin{equation*}
x_1=\phi^l(x_1)=x_{l+1}
\end{equation*}
holds. This is a contradiction.  \begin{flushright}  $\Box$  \end{flushright}

\section{Proof of Theorem 1}
Before stating the proof of Theorem 1, we introduce the notion of contact Hamiltonian function. Let ${M}$ be a smooth manifold without boundary and let ${\alpha\in \Omega^1(M)}$ be a contact form on $M$ (${\xi=\ker(\alpha)}$). A Reeb vector field ${R\in \Gamma(TM)}$ is defined as follows.

\begin{gather*}
\alpha(R)=1  \\
d\alpha(R,\cdot)=0 
\end{gather*}

For any smooth function ${h\in C_c^{\infty}(M)}$, there exists an unique contact vector field ${X_h\in \Gamma_{\xi}^c(TM)}$ which satisfies the following condition.

\begin{gather*}
X_h=h\cdot R+Z \ \ (Z\in \xi)
\end{gather*}

In fact, $X_h$ is contact vector field if and only if ${\mathcal{L}_{X_h}(\alpha)|_{\xi}=0}$ holds ($\mathcal{L}$ is the Lie derivative.). So,
\begin{equation*}
\mathcal{L}_{X_h}(\alpha)(Y)=dh(Y)+d\alpha(X_h,Y)=dh(Y)+d\alpha(Z,Y)=0
\end{equation*}

holds for any ${Y\in \xi}$. Because ${d\alpha}$ is non-degenerate on ${\xi}$, above equation determines ${Z\in \xi}$ uniquely. ${X_h}$ is a contact vector field associated to a contact Hamiltonian function $h$. We denote the time $t$ flow of ${X_h}$ by ${\phi^t_h}$ and time $1$ flow of ${X_h}$ by ${\phi_h}$.

Let ${(M,\xi)}$ be a contact manifold without boundary. We fix a point ${p\in (M,\xi)}$ and a sufficiently small open neighborhood ${U\subset M}$ of $p$. Let ${(x_1,y_1,\cdots,x_n,y_n,z)}$ be a coordinate of ${\mathbb{R}^{2n+1}}$. Let ${\alpha_0\in \Omega^1(\mathbb{R}^{2n+1})}$ be the following contact form on ${\mathbb{R}^{2n+1}}$.
\begin{equation*}
\alpha_0=\frac{1}{2}\sum_{1\le i\le n}(x_idy_i-y_idx_i)+dz
\end{equation*}

By using famous Moser's arguments, we can assume that there exists an open neighborhood of the origin ${V\subset \mathbb{R}^{2n+1}}$ and a diffeomorphism
\begin{gather*}
F:V\longrightarrow U
\end{gather*}
which satisfies the following condition.
\begin{equation*}
\xi |_U=\textrm{Ker}(F_*\alpha_0)
\end{equation*}
So, we first prove the theorem for ${(M,\xi)=(V,\ker(\alpha_0))}$ and apply this to ${(M,\xi)}$. 

We fix ${k\in \mathbb{N}_{\ge 1}}$ and ${R>0}$ so that 
\begin{equation*}
\{(x_1,y_1,\cdots,z)\in \mathbb{R}^{2n+1} \ | \ |(x_1,\cdots,y_n)|<R, |z|<R\}\subset V
\end{equation*}
holds. Let ${f\in C^{\infty}_c(V)}$ be a contact Hamiltonian function. Then its contact Hamiltonian vector field ${X_f}$ can be written in the following form.
\begin{gather*}
X_f(x_1,\cdots,z)=\sum_{1\le i\le n}(-\frac{\partial f}{\partial y_i}+\frac{x_i}{2}\frac{\partial f}{\partial z})\frac{\partial}{\partial x_i}   \\
+\sum_{1\le i\le n}(\frac{\partial f}{\partial x_i}+\frac{y_i}{2}\frac{\partial f}{\partial z})\frac{\partial}{\partial y_i}  \\
+(f-\sum_{1\le i\le n}\frac{x_i}{2}\frac{\partial f}{\partial x_i}-\sum_{1\le i\le n}\frac{y_i}{2}\frac{\partial f}{\partial y_i})\frac{\partial}{\partial z}
\end{gather*}

Let ${e:\mathbb{R}^{2n}\longrightarrow \mathbb{R}}$ be a quadric function as follows.
\begin{equation*}
e(x_1,y_1,\cdots,x_n,y_n)=x_1^2+y_1^2+\sum_{2\le i\le n}\frac{x_i^2+y_i^2}{2}
\end{equation*}
We consider the following contact Hamiltonian function on ${\mathbb{R}^{2n+1}}$.
\begin{equation*}
h(x_1,y_1,\cdots,x_n,y_n,z)=\beta(z)\rho(e(x_1,y_1,\cdots,x_n,y_n))
\end{equation*}
Here, ${\beta:\mathbb{R}\rightarrow [0,1]}$ and ${\rho:\mathbb{R}_{\ge 0}\rightarrow \mathbb{R}_{\ge 0}}$ are smooth functions which satisfy the following conditions.
\begin{itemize}
\item ${\textrm{supp}(\rho)\subset [0,\frac{R^2}{2}]}$
\item ${\rho(r)\ge \rho'(r)\cdot r}$, ${-\frac{\pi}{2k}<\rho'(r)\le\frac{\pi}{2k}}$
\item There exists an unique ${a\in[0,\frac{R^2}{2}]}$ which satisfies the following conditions.
\begin{equation*}
\begin{cases}  \rho'(r)=\frac{\pi}{2k}\Longleftrightarrow r=a \\
\rho(a)=\frac{\pi}{2k}\cdot a
\end{cases}
\end{equation*}
\item $\textrm{supp}(\beta)\subset [-\frac{R}{2},\frac{R}{2}]$
\item $\beta(0)=1$, $\beta'(0)=0$, ${\beta^{-1}(1)=0}$
\end{itemize}

Then, we can prove the following lemma.
\begin{Lem}
Let ${h\in C^{\infty}_c(V)}$ be a contact Hamiltonian function as above. Then, 
\begin{equation*}
(q,\phi_h(q),\cdots,\phi_h^{2k-1}(q))\in P^{2k}(\phi_h)
\end{equation*}
holds if and only if 
\begin{equation*}
q\in \{(x_1,y_1,0,\cdots,0)\in V \ | \ x_1^2+y_1^2=a\}=S_a
\end{equation*}
holds.
\end{Lem}

In order to prove this lemma, we first calculate the behavior of the function ${z(\phi^t_h(q))}$ for fixed ${q\in V}$ (Here, $z$ is the ${(2n+1)}$-th coordinate of ${\mathbb{R}^{2n+1}}$.).
\begin{gather*}
\frac{d}{dt}(z(\phi^t_h(q)))=h-\sum_{1\le i\le n}\frac{x_i}{2}\frac{\partial h}{\partial x_i}-\sum_{1\le i\le n}\frac{y_i}{2}\frac{\partial h}{\partial y_i} \\
=\beta(z)\{\rho(e)-\sum_{1\le i\le n}\frac{x_i}{2}\frac{\partial}{\partial x_i}(\rho(e))-\sum_{1\le i\le n}\frac{y_i}{2}\frac{\partial}{\partial y_i}(\rho(e)) \} \\
=\beta(z)\{\rho(e)-\rho'(e)\cdot e\} \ge 0
\end{gather*}

So, this inequality implies that 
\begin{equation*}
\phi_h^{2k}(q)=q \Longrightarrow  \frac{d}{dt}(z(\phi^t_h(q)))=0
\end{equation*}
holds. 

Next, we study the behavior of ${x_i(\phi^t_h(q))}$ and ${y_i(\phi^t_h(q))}$. Let ${\phi_i}$ be the following projection.
\begin{gather*}
\pi_i:\mathbb{R}^{2n+1}\longrightarrow \mathbb{R}^2  \\
(x_1y_1,\cdots,x_n,y_n,z)\mapsto (x_i,y_i)
\end{gather*}
Then, ${Y_h^i=\pi_i(X_h)}$ can be decomposed into angular component ${Y_h^{i,\theta}}$ and radius component ${Y_h^{i,r}}$ as follows. 
\begin{gather*}
Y_h^{i,\theta}(x_1,y_1,\cdots,z)=-\frac{\partial h}{\partial y_i}\frac{\partial}{\partial x_i}+\frac{\partial h}{\partial x_i}\frac{\partial }{\partial y_i} \\
Y_h^{i,r}(x_1,y_1,\cdots,z)=(\frac{1}{2}\frac{\partial h}{\partial z})(x_i\frac{\partial}{\partial x_i}+y_i\frac{\partial}{\partial y_i})
\end{gather*}
Let ${w_i}$ be the complex coordinate of ${(x_i,y_i)}$ (${w_i=x_i+\sqrt{-1}y_i}$). Then, the angular component causes the following rotation on ${w_i}$ (If we ignore the $z$-coordinate).
\begin{equation*}
\arg (w_i)\longrightarrow \arg(w_i)+2\rho'(e(x_1,\cdots,y_n))\beta(z)C_it
\end{equation*}
\begin{equation*}
C_i=\begin{cases} 1 & i=1 \\ \frac{1}{2} & 2\le i\le n   \end{cases}
\end{equation*}
By the assumptions of ${\rho}$ and ${\beta}$, ${|2\rho'(e(x_1,\cdots,y_n))\beta(z)C_i|}$ is at most ${\frac{2\pi}{2k}}$ and equality holds if and only if ${(x_1,y_1,\cdots,x_n,y_n,z)\in S_a}$ holds. On the circle ${q\in S_a}$, ${\phi_h}$ is the ${\frac{2\pi}{2k}}$-rotation of the circle ${S_a}$. This implies that the lemma holds. \begin{flushright}  $\Box$  \end {flushright}

Next, we perturb the contactomorphism ${\phi_h}$. Let ${(r,\theta)}$ be a coordinate of ${(x_1,y_1)\in \mathbb{R}^2\backslash (0,0)}$ as follows.
\begin{equation*}
x_1=r\cos \theta, \ \ y_1=r\sin \theta
\end{equation*}

We fix ${\epsilon>0}$. Then ${\epsilon(1-\cos(k\theta))}$ is a contact Hamiltonian function on ${\mathbb{R}^2\backslash (0,0)\times \mathbb{R}^{2n-1}}$ and its contact Hamiltonian vector field can be written in the following form.
\begin{gather*}
X_{\epsilon(1-\cos (k\theta))}=-\frac{\epsilon k}{r}\sin (k\theta)\frac{\partial}{\partial r}+\epsilon(1-\cos(k\theta))\frac{\partial}{\partial z}
\end{gather*}

So ${\phi_{\epsilon(1-\cos(k\theta))}}$ only changes the $r$ of ${(x_1,y_1)}$-coordinate and ${z}$-coordinate as follows.
\begin{equation*}
(r,\theta,x_2,y_2,\cdots,x_n,y_n,z)\mapsto (\sqrt{r^2-2\epsilon k\sin(k\theta)},\theta,x_2,\cdots,y_n,z+\epsilon(1-\cos(k\theta)))
\end{equation*}

We fix two small open neighborhood of the circle ${S_a}$ as follows.
\begin{gather*}
S_a\subset W_1\subset W_2\subset \mathbb{R}^2\backslash (0,0)\times \mathbb{R}^{2n-1}
\end{gather*}
We also fix a cut-off function ${\eta:\mathbb{R}^{2n+1}\rightarrow [0,1]}$ as follows.
\begin{gather*}
\eta((x_1,\cdots,z))=1 \ \  \ \ \ \ ((x_1,\cdots,z)\in W_1) \\
\eta((x_1,\cdots,z))=0 \ \  \ \ \ \  ((x_1,\cdots,z)\in \mathbb{R}^{2n+1}\backslash W_2)
\end{gather*}
Then, ${\eta(x_1.\cdots,z)\cdot \epsilon(1-\cos(k\theta))}$ is defined on ${\mathbb{R}^{2n+1}}$. We denote this contact Hamiltonian function by ${g_{\epsilon}}$. We define ${\phi_{\epsilon}\in \textrm{Cont}^c_0(\mathbb{R}^{2n+1},\ker(\alpha_0))}$ by the composition ${\phi_{g_{\epsilon}}\circ \phi_h}$. 

\begin{Lem}
We take ${\epsilon>0}$ sufficiently small. We define ${2k}$ points ${\{a_i\}_{1\le i\le 2k}}$ as follows.
\begin{equation*}
a_i=(\sqrt{a}\cos (\frac{i\pi}{k}),\sqrt{a}\sin (\frac{i\pi}{k}),0,\cdots,0))\in S_a
\end{equation*}
Then ${P^{2k}(\phi_{\epsilon})}$ has only one point ${[a_1,a_2,\cdots,a_{2k}]}$.
\end{Lem}

The proof of this lemma is as follows. On $W_1$, $\phi_{g_{\epsilon}}$ only changes the $r$-coordinate of ${(x_1,y_1)}$ and ${z}$-coordinate. On the circle ${S_a}$, the fixed points of ${\phi_{g_{\epsilon}}}$ are 2k points ${\{a_i\}}$. This implies that on $W_1$, ${\phi_{\epsilon}}$
increases the angle ${\theta}$ at most ${\frac{2\pi}{2k}}$ and the equality holds on only 2k points ${\{a_i\}}$. From the arguments in the proof of Lemma 1, this implies that 
\begin{equation*}
P^{2k}(\phi_{\epsilon})=\{[a_1,a_2,\cdots,a_{2k}]\}
\end{equation*}
holds if $\epsilon$ is sufficiently small.  \begin{flushright}  $\Box$  \end{flushright}

Finally, we prove Theorem 1. We define ${\psi_{\epsilon}\in \textrm{Cont}^c_0(M,\xi)}$ as follows.

\begin{equation*}
\psi_{\epsilon}(x)=\begin{cases} F\circ \phi_{\epsilon}\circ F^{-1}(x)  &  x\in U   \\ x & x\in M\backslash U\end{cases}
\end{equation*}

Lemma 2 implies that 
\begin{equation*}
P^{2k}(\psi_{\epsilon})=\{[F(a_1),\cdots,F(a_{2k})]\}
\end{equation*}
holds. Proposition 1 implies that ${\psi_{\epsilon}\in \textrm{Cont}^c_0(M,\xi)\backslash \textrm{Aut}(M,\xi)}$ holds. Because ${p\in M}$ is any point and ${U}$ is any small open neighborhood of $p$, we proved Theorem 1.  \begin{flushright}  $\Box$ \end{flushright}

\section{Proof of Theorem 2}
Let $M$ be a $m$-dimensional smooth manifold without boundary. We fix a point ${p\in M}$. Let ${U}$ be an open neighborhood of $p$ and let ${V\subset \mathbb{R}^m}$ be an open neighborhood of the origin such that there is a diffeomorphism 
\begin{equation*}
F:V\longrightarrow U .
\end{equation*}
In order to prove Theorem 2, it suffices to prove that there exists a sequence ${\psi_{\epsilon}}$ so that 
\begin{itemize}
\item ${\psi_{\epsilon}\in \textrm{Diff}^c_0(M)\backslash \textrm{Aut}(M)}$
\item ${\textrm{supp}(\Psi_{\epsilon})\subset U}$
\item ${\psi_{\epsilon}\longrightarrow} \textrm{Id}$ as ${\epsilon \rightarrow 0}$
\end{itemize}
hold. 

First, assume that ${m}$ is odd (${m=2n+1}$). In this case, $\alpha_0$ is a contact form on $V$. Let ${\phi_{\epsilon}}$ be a contactomorphism which we constructed in the proof of Theorem 1.
\begin{itemize}
\item ${\phi_{\epsilon}\in \textrm{Cont}^c_o(\mathbb{R}^{2n+1},\ker(\alpha_0))}$
\item ${\textrm{supp}(\phi_{\epsilon})\subset V}$
\end{itemize}
We define ${\psi_{\epsilon}\in \textrm{Diff}^c_0(M)}$ as follows.
\begin{equation*}
\psi_{\epsilon}(x)=\begin{cases} F\circ \phi_{\epsilon}\circ F^{-1}(x)   & x\in U \\ x & x\in M\backslash U \end{cases}
\end{equation*}

Then, ${\sharp P^{2k}(\psi_{\epsilon})=1}$ holds and this implies that ${\psi_{\epsilon}\in \textrm{Diff}^c_0(M)\backslash \textrm{Aut}(M)}$ holds.

Next, assume that $m$ is even (${m=2n}$). Let ${\omega_0}$ be a standard symplectic form on ${(x_1,y_1,\cdots,x_n,y_n)\in \mathbb{R}^{2n}}$ as follows.
\begin{equation*}
\omega_0=\sum_{1\le i\le n}dx_i\wedge dy_i
\end{equation*}
As in the contact case, Hamiltonian function ${h\in C_c^{\infty}(V)}$ defines Hamiltonian vector field ${X_h}$ as follows.  
\begin{equation*}
\omega_0(X_h,\cdot)=-dh
\end{equation*}
We denote the time $1$ flow of $X_h$ by ${\phi_h}$. 

\begin{Rem}
We can use the arguments in \cite{AF}. However, our perturbation is slightly different from that of \cite{AF} (Our perturbation is in the relatively compact domain U.). So, we explain our arguments in detail for the sake of self-containedness.
\end{Rem}

We assume that ${B(R)\subset V}$ holds (${B(R)}$ is a radius $R$-ball in ${\mathbb{R}^{2n}}$.). Let ${e:\mathbb{R}^{2n}\rightarrow \mathbb{R}}$ be a following quadric function.
\begin{equation*}
e(x_1,y_1,\cdots,x_n,y_n)=x_1^2+y_1^2+\sum_{2\le i\le n}\frac{x_i^2+y_i^2}{2}
\end{equation*}
Let ${\rho:[0,\infty)\rightarrow \mathbb{R}_{\ge 0}}$ be a non-negative function which satisfies the following conditions.
\begin{itemize}
\item ${\textrm{supp}(\rho)=[0,\frac{R^2}{2}]}$
\item ${-\frac{k}{2\pi}<\rho'\le \frac{k}{2\pi}}$, $\rho'(r)=\frac{k}{2\pi}\Longleftrightarrow r=a$
\end{itemize}
Let ${h\in C^{\infty}_c(V)}$ be a Hamiltonian function as follows.
\begin{equation*}
h(x_1,y_1,\cdots,x_n,y_n)=\rho(e(x_1,y_1,\cdots,x_n,y_n))
\end{equation*}
The following lemma is an analogy of Lemma 1. 
\begin{Lem}
Let $h\in C^{\infty}_c(V)$ be a Hamiltonian function as above. Then, 
\begin{equation*}
(q,\phi_h(q),\cdots,\phi_h^{2k-1}(q))\in P^{2k}(\phi_h)
\end{equation*}
holds if and only if
\begin{equation*}
q\in \{(x_1,y_1,0\cdots,0)\in V \ | \ x_1^2+y_1^2=a \}=S_a
\end{equation*}
holds.
\end{Lem}
The proof of this lemma is similar to that of Lemma 1. Let ${w_i}$ be a complex coordinate of ${(x_i,y_i)}$ (${w_i=x_i+\sqrt{-1}y_i}$). Then, ${X_h}$ causes the following rotaion on ${w_i}$.
\begin{equation*}
w_i\mapsto \exp(2\rho'(e(x_1,\cdots,y_n))C_it\sqrt{-1})\cdot w_i
\end{equation*}
\begin{equation*}
C_i=\begin{cases} 1 & i=1 \\ \frac{1}{2}  & 2\le i\le n \end{cases}
\end{equation*}
${|2\rho'(e)C_i|}$ is at most ${\frac{2\pi}{2k}}$ and equality holds if and only if ${(x_1,\cdots,y_n)\in S_a}$ holds. This implies that the lemma holds. \begin{flushright}   ${\Box}$ \end{flushright}

Next we perturb this ${\phi_h\in \textrm{Ham}^c(V)}$. Let ${(r,\theta)}$ be a coordinate of ${\mathbb{R}^2\backslash (0,0)}$ as follows.
\begin{equation*}
x_1=r\cos \theta, y_1=r\sin \theta
\end{equation*}
We fix ${\epsilon >0}$, then ${\epsilon \cos(k\theta)}$ is a Hamiltonian function on ${\mathbb{R}^2\backslash (0,0)\times \mathbb{R}^{2n-2}}$. The Hamiltonian vector field of this Hamiltonian function can be written in the following form.
\begin{equation*}
X_{\epsilon \cos(k\theta)}=\frac{\epsilon k}{r}\sin(k\theta)\frac{\partial}{\partial r}
\end{equation*}
So, ${\phi_{\epsilon \cos(k\theta)}}$ only changes the $r$ of ${(x_1,y_1)}$-coordinate. 
\begin{equation*}
(r,\theta,x_2,y_2,\cdots,x_n,y_n)\mapsto (\sqrt{r^2+2\epsilon k\sin(k\theta)},\theta,x_2,y_2,\cdots,x_n,y_n)
\end{equation*}
We fix two small open neighborhoods of ${S_a}$ and cut-off function ${\eta:\mathbb{R}^{2n}\rightarrow \mathbb{R}}$ as follows.
\begin{gather*}
S_a\subset W_1\subset W_2\subset \mathbb{R}^2\backslash (0,0)\times \mathbb{R}^{2n-2} \\
\eta(x)=\begin{cases} 1  & x\in W_1  \\ 0 & x\in \mathbb{R}^{2n}\backslash W_2  \end{cases}
\end{gather*}
We denote the Hamiltonian function ${\eta \cdot \epsilon \cos(k\theta)\in C^{\infty}_c(\mathbb{R}^{2n})}$ by ${g_{\epsilon}}$ and let ${\phi_{\epsilon}\in \textrm{Ham}^c(\mathbb{R}^{2n})}$ be the composition ${\phi_{g_{\epsilon}}\circ \phi_h}$. The next lemma is an analogy of Lemma 2 and its proof is almost the same. So we omit to prove it.
\begin{Lem}
We take ${\epsilon>0}$ sufficiently small. We define ${2k}$ points ${\{a_i\}_{1\le i\le 2k}}$ as follows.
\begin{equation*}
a_i=(\sqrt{a}\cos(\frac{i\pi}{k}),\sqrt{a}\sin(\frac{i\pi}{k}),0,\cdots,0)\in S_a
\end{equation*}
Then, ${P^{2k}(\phi_{\epsilon})}$ has only one point ${[a_1,a_2,\cdots,a_{2k}]}$.
\end{Lem}

We define ${\psi_{\epsilon}\in \textrm{Diff}^c_0(M)}$ as follows.
\begin{equation*}
\psi_{\epsilon}=\begin{cases} F\circ \phi_{\epsilon}\circ F^{-1}   &  x\in U  \\  x  &  x\in M\backslash U  \end{cases}
\end{equation*}
Lemma 4 implies that ${\sharp P^{2k}(\psi_{\epsilon})=1}$ holds and this implies that
\begin{equation*}
\psi_{\epsilon}\in \textrm{Diff}^c_0(M)\backslash \textrm{Aut}(M)
\end{equation*}
holds. So we proved Theorem 2. \begin{flushright}   ${\Box}$  \end{flushright}


\begin{thebibliography}{9}
\bibitem{AF}
P. Albers, U. Frauenfelder.
\textit{Square roots of Hamiltonian diffeomorphisms.}
J. Symplectic. Geom. Volume 12, Number 3 (2014), 427-434
\bibitem{M}
J. Milnor.
\textit{Remarks on infinite-dimensional Lie groups}.
relativity, groups and topology Ⅱ, COURSE 10
\bibitem{PS}
L. Polterovich, E. Schelukhin.
\textit{Autonomous Hamiltonian flows, Hofer's geometry and persistence modules.}
Selecta Mathematica January 2016, Volume 22, Issue 1, pp227-296
\bibitem{S}
Y. Sugimoto.
\textit{Spectral spread and non-autonomous Hamiltonian diffeomorphisms.}
manuscripta math. (2018). https://doi.org/10.1007/s00229-018-1078-0 
\end{thebibliography}
\end{document}